\documentclass[12pt]{article}
\usepackage{amsmath,amssymb,amsthm}
\usepackage{mathrsfs}
\usepackage{graphicx}
\usepackage{multirow}
\usepackage[utf8]{inputenc}
\usepackage[english]{babel}
\usepackage{geometry}
\usepackage{cite}
\usepackage{hyperref}
\usepackage[pagewise]{lineno}
\usepackage{enumitem}
\usepackage{bm}

\newtheorem{theorem}{Theorem}[section]
\newtheorem{lemma}[theorem]{Lemma}

\newtheorem{corollary}[theorem]{Corollary}

\newtheorem{remark}[theorem]{Remark}

\newcommand{\vol}{\text{vol}}

\title{Gradient Estimates for Diffusion Equations on Weighted Graphs  }
\author{
Shoudong Man\thanks{College of Science and Technology, Tianjin University of Finance and Economics, Tianjin 300222,
China. Corresponding author:  Email: manshoudong@163.com; shoudongmantj@tjufe.edu.cn.}
}

\date{}

\begin{document}

\maketitle

\begin{abstract}
We establish a unified approach to Li–Yau-type gradient estimates for diffusion equations on graphs, incorporating both fractional and standard Laplacians with variable exponents. Based on novel discrete chain rules, our estimates yield Harnack inequalities,
from which we derive heat kernel bounds, volume growth estimates and Buser-type isoperimetric inequalities. This work provides a foundational framework for analyzing nonlocal diffusion on discrete structures.
\end{abstract}



\section{Introduction}

The celebrated Li-Yau inequality, a cornerstone of geometric analysis established in \cite{LI-YAU}, provides a powerful gradient estimate for positive solutions to the heat equation on Riemannian manifolds. In recent years, there has been significant interest in extending these results to discrete settings, particularly on graphs, such as those in \cite{Bauer2015, Anal,Lin2017, Lin1}. While significant progress has been made for linear equations, the study of nonlinear diffusion equations with variable exponents, especially those involving fractional Laplace operators, remains largely undeveloped. This paper aims to bridge this gap by establishing a unified framework for gradient estimates on graphs.

In 2025, M. Zhang, Y. Lin and Y. Yang investigated fractional diffusion equations on graphs \cite{ZhangLinYang1,ZhangLinYang2,ZhangLinYang3}.
They developed a comprehensive groundwork  for fractional Sobolev spaces and fractional Laplace on graphs, and established the existence of solutions to fractional Laplace equations on graphs. These foundational existence results provide a crucial starting point for further analysis of nonlinear fractional equations on graphs.

In this paper, we mainly investigate the general fractional diffusion equations on graphs of the form
\begin{align}
\Delta^{\epsilon} (u^{m(x,t)}) - b(x,u)u = \delta(x,t) u_t + \gamma(t) u^{m(x,t)}, \label{main02}
\end{align}
where $\Delta^{\epsilon}$ denotes the fractional Laplace operator on graphs with $\epsilon\in(0,1)$,
$m(x,t)$ is a variable exponent capturing spatially and temporally heterogeneous diffusion properties, $b(x,u)$ represents a convection term, $\delta(x,t)$ is a time-dependent coefficient, and $\gamma(t)$ is a given function.
We also consider the corresponding equation with the standard graph Laplacian
\begin{align}
\Delta (u^{m(x,t)}) - b(x,u)u = \delta(x,t) u_t + \gamma(t) u^{m(x,t)}. \label{main01}
\end{align}
We develop a discrete chain rule for variable exponents on graphs, which serves as a fundamental tool for handling nonlinear and nonlocal terms.
As a result, we establish unified global gradient estimates for positive solutions to both fractional and standard diffusion equations \eqref{main02} and \eqref{main01}. From these gradient estimates, we derive Harnack inequalities, which in turn yield heat kernel estimates, volume growth bounds and Buser-type inequalities.

The paper is organized as follows:
Section 2 provides the preliminaries and notation.
Section 3 establishes discrete chain rules and the main gradient estimates.
Section 4 derives Harnack inequalities from the estimates.
Section 5 presents applications of the results.

\section{Preliminaries and Notation}
Following \cite{ZhangLinYang1}, we define the fractional framework on graphs.
Let \(G=(V,E,\vartheta,w)\) denote a connected, locally finite and stochastically complete graph equipped with vertex set \(V\), edge set \(E\), measure \(\vartheta\) and weight \(w\). A graph \(G\) is called stochastically complete if its heat kernel \(p:[0,+\infty)\times V\times V\rightarrow\mathbb{R}\) satisfies
\[
\sum_{y\in V}p_{t}(x,y)\vartheta(y)=1,\quad\forall t>0,\;\forall x\in V.
\]
Conversely, \(G\) is termed stochastically incomplete if there exist some \(t_{0}>0\) and \(x_{0}\in V\)
such that \(\sum_{y\in V}p_{t_{0}}(x_{0},y)\vartheta(y)<1\). Denote $x\sim y$ if vertex $x$ is adjacent to vertex $y$.
Let $w_{xy}$ denote an edge weight satisfying $w_{xy} = w_{yx} > 0$ for $x \sim y$ and $w_{xy} = 0$ otherwise.
The degree of vertex $x$ is denoted by $deg(x)=\sum_{y\sim x}\omega_{xy}.$
A graph is called locally finite if each vertex has finite neighbors.
A graph $G$ is called connected if for any vertices $x,y\in V$, there exists a sequence $\{x_{i}\}_{i=0}^{n}$ satisfying
$$x=x_{0}\sim x_{1}\sim x_{2}\sim \cdot\cdot\cdot \sim x_{n}=y.$$
Throughout this paper, we assume that all graphs are connected.
The distance $d(x,y)$ between two points $x$ and $y$ is defined as the length of the shortest path connecting them, that is,
\begin{align*}
d(x,y)=\min\{\eta: ~ there ~exists ~a ~path ~x=x_{0}\sim x_{1}\sim x_{2}\sim \cdot\cdot\cdot \sim x_{\eta}=y\}.
\end{align*}

 Assume $\epsilon\in(0,1)$. Define a kernel as
\begin{equation}
W_{\epsilon}(x,y)=\frac{\epsilon}{\Gamma(1-\epsilon)}\vartheta(x)\vartheta(y)\int_{0}^{+\infty}p_{t}(x,y)t^{-1-\epsilon}dt,\quad \forall x\neq y\in V, \notag
\end{equation}
where $\Gamma(\cdot)$ denotes the Gamma function and $p_{t}(x,y)$ is the heat kernel. Moreover, from \cite{ZhangLinYang1}
$$\sum_{y\in V, y \neq x} W_{\epsilon}(x,y)\leq C_{x,\epsilon},$$
where $C_{x,\epsilon}$ is a positive constant depending only on $x$ and $\epsilon$.

From \cite{ZhangLinYang1}, for any function $u:V\rightarrow \mathbb{R}$, the fractional Laplace of $u$ can be expressed as
\begin{align}
\Delta^{\epsilon} u(x) = \frac{1}{\vartheta(x)} \sum_{y\in V, y \neq x} W_{\epsilon}(x,y) \left(u(y) - u(x)\right). \label{fenshu}
\end{align}
The associated gradient form reads
 \begin{align}
\Gamma^{\epsilon}(u,v)(x)=\frac{1}{2\vartheta(x)}\sum_{y\in V, y \neq x} W_{\epsilon}(x,y)(u(y)-u(x))(v(y)-v(x)).\notag
\end{align}
We denote $\Gamma^{\epsilon}(u)=\Gamma^{\epsilon}(u,u)$ for short.
The length of the gradient for $u$ is denoted by
\begin{equation}
|\nabla^{\epsilon}u|(x)=\sqrt{\Gamma^{\epsilon}(u)(x)}=\left(\frac{1}{2\vartheta(x)}\sum_{y\in V,y\neq x}W_{\epsilon}(x,y)(u(x)-u(y))^{2}\right)^{\frac{1}{2}}. \notag
\end{equation}

From \cite{YLS}, the standard $\vartheta$-Laplacian of $u$ is defined as
 \begin{align*}
\Delta u(x)=\frac{1}{\vartheta(x)}\sum_{y\sim x}\omega_{xy}\left(u(y)-u(x)\right).
\end{align*}
 The associated gradient form reads
 \begin{align}
\Gamma(u,v)(x)=\frac{1}{2\vartheta(x)}\sum_{y\sim x}\omega_{xy}(u(y)-u(x))(v(y)-v(x)).\notag
\end{align}
We denote $\Gamma(u)=\Gamma(u,u)$ for short.
The length of the gradient for $u$ is denoted by
\begin{equation*}
|\nabla u|(x)=\sqrt{\Gamma(u)(x)}=\Big{(}\frac{1}{2\vartheta(x)}\sum_{y\sim x} \omega_{xy}\big{(}u(y)-u(x)\big{)}^{2}\Big{)}^{1/2}.
\end{equation*}
For the Laplacian, we assume the following bounds
\begin{align*}
&\omega_{\min}:=\inf_{x,y\in V} W_{\epsilon}(x,y)<\infty,
D^{\epsilon}_{\vartheta}:=\max_{x\in V}\frac{\sum_{y\in V,y\neq x}W_{\epsilon}(x,y)}{\vartheta(x)}<\infty,  \\
&D_{\vartheta}:=\max_{x\in V}\frac{\deg(x)}{\vartheta(x)}<\infty,
\vartheta_{\max}:=\sup_{x\in V}\vartheta(x)<\infty.
\end{align*}

\section{Gradient Estimates}

In this section, we derive gradient estimates for positive solutions to equations \eqref{main02} and \eqref{main01}. Our approach begins with establishing discrete chain rule identities for variable exponents, which are fundamental to our analysis.

\subsection{Discrete Chain Rules}

\begin{lemma}[Discrete Chain Rule for the Fractional Laplacian] \label{lem:chain_rule1}
Let $G = (V,E,\vartheta,w)$ be a connected, locally finite, and stochastically complete graph equipped
with the fractional Laplace operator $\Delta^{\epsilon}$. For any positive function $u: V \to \mathbb{R}^+$
and any variable exponent $m: V \times [0,\infty) \to \mathbb{R}^+$, the following identity holds at each vertex $x \in V$
\[
\Delta^{\epsilon} (u^{m(x,t)}) = 2 u^{\frac{m(x,t)}{2}} \Delta^{\epsilon} (u^{\frac{m(x,t)}{2}}) + 2 \Gamma^{\epsilon}(u^{\frac{m(x,t)}{2}}).
\]
\end{lemma}

\begin{proof}
 For any vertex $x \in V$, we have
\begin{align*}
\Delta^{\epsilon} (u^{m(x,t)})(x) &= \frac{1}{\vartheta(x)} \sum_{y\in V, y \neq x} W_{\epsilon}(x,y) (u^{m(y,t)}(y) - u^{m(x,t)}(x)) \\
&= \frac{1}{\vartheta(x)} \sum_{y\in V, y \neq x} W_{\epsilon}(x,y) \left[ \left(u^{\frac{m(y,t)}{2}}(y)\right)^2 - \left(u^{\frac{m(x,t)}{2}}(x)\right)^2 \right].
\end{align*}
Using the identity $a^2 - b^2 = 2b(a - b) + (a - b)^2$, we get
\begin{align*}
\Delta^{\epsilon} (u^{m(x,t)})(x) &= \frac{2}{\vartheta(x)} \sum_{y\in V, y \neq x} W_{\epsilon}(x,y) u^{\frac{m(x,t)}{2}}(x)
\left(u^{\frac{m(y,t)}{2}}(y) - u^{\frac{m(x,t)}{2}}(x)\right) \\
&\quad + \frac{1}{\vartheta(x)} \sum_{y\in V, y \neq x} W_{\epsilon}(x,y) \left(u^{\frac{m(y,t)}{2}}(y) - u^{\frac{m(x,t)}{2}}(x)\right)^2 \\
&= 2 u^{\frac{m(x,t)}{2}}(x) \Delta^{\epsilon} (u^{\frac{m(x,t)}{2}})(x) + 2 \Gamma^{\epsilon}(u^{\frac{m(x,t)}{2}})(x).
\end{align*}
This completes the proof.
\end{proof}

The chain rule for the standard Laplacian $\Delta$ can be derived similarly.

\begin{lemma}[Discrete Chain Rule for the Standard Laplacian] \label{lem:chain_rule}
Let $G = (V,E)$ be a finite or locally finite graph equipped with the standard Laplace operator $\Delta$.
For any positive function $u: V \to \mathbb{R}^+$ and any variable exponent $m: V \times [0,\infty) \to \mathbb{R}^+$,
the following identity holds at each vertex $x \in V$
\[
\Delta (u^{m(x,t)})(x) = 2u^{\frac{m(x,t)}{2}}(x)\Delta (u^{\frac{m(x,t)}{2}})(x) + 2\Gamma(u^{\frac{m(x,t)}{2}})(x).
\]
\end{lemma}

\begin{remark}
The discrete chain rule identities established in Lemmas $\ref{lem:chain_rule1}$ and $\ref{lem:chain_rule}$ reveal fundamental differences between graph and classical analysis. These differences are not merely technical artifacts but reflect deep geometric distinctions between discrete graph structures and continuous manifolds.

The discrete chain rule identities serve as a cornerstone in extending the Li-Yau type theory to graphs and enable the development of Harnack inequalities and fundamental solution estimates in this discrete setting. These identities initiate our derivation of gradient estimates for positive solutions to equations $\eqref{main02}$ and $\eqref{main01}$.
\end{remark}

\subsection{ Gradient Estimates}

We now establish our main gradient estimate for positive solutions to equation \eqref{main02}.
\begin{theorem}[Gradient Estimate for the fractional Laplacian in the Variable Exponents case]\label{thm:gradient_linear}
Let $G = (V,E,\vartheta,w)$ be a connected, locally finite, and stochastically complete graph. Suppose $u = u(x,t) > 0$ is a solution to equation \eqref{main02}, where
$m(x,t): V \times (0, \infty) \to \mathbb{R}^+$ is a variable exponent,
  $\delta(x,t): V \times (0, \infty) \to \mathbb{R}$ is a weight function,
  $b(x,u) \leq \varphi(x,t) u^{m(x,t)-1}$ for some function $\varphi(x,t): V \times (0, \infty) \to \mathbb{R}$, and
$\gamma(t) \leq c_1 + \dfrac{c_2}{t}$ for constants $c_1, c_2 \in \mathbb{R}$.
Then for any $x \in V$, $t > 0$, and any function $\alpha(t) \geq 1$, the following gradient estimate holds
\begin{align*}
\frac{2\Gamma^{\epsilon}(u^{\frac{m(x,t)}{2}})}{u^{m(x,t)}} - \alpha(t)\frac{\delta u_t}{u^{m(x,t)}}
\leq 2\alpha(t)D^{\epsilon}_{\vartheta} + \alpha(t)\varphi(x,t) + \alpha(t)\left(c_1 + \frac{c_2}{t}\right).
\end{align*}
\end{theorem}

\begin{proof}
Define the function
\[
F(x,t) = \frac{2\Gamma^{\epsilon}(u^{\frac{m(x,t)}{2}})}{u^{m(x,t)}} - \alpha(t)\cdot\frac{\delta u_t}{u^{m(x,t)}}.
\]
From equation \eqref{main02}, we have
\begin{align}
\delta u_t=\Delta^{\epsilon} (u^{m(x,t)}) - b(x,u)u -\gamma(t)u^{m(x,t)}.\label{TC2}
\end{align}
By \eqref{TC2} and Lemma \ref{lem:chain_rule1}, we obtain
\begin{align}
F(x,t) = & \frac{2\Gamma^{\epsilon}(u^{\frac{m(x,t)}{2}})}{u^{m(x,t)}} - \alpha(t)\frac{\delta u_t}{u^{m(x,t)}} \notag\\
=& \frac{2\Gamma^{\epsilon}(u^{\frac{m(x,t)}{2}})-\alpha(t)\Delta^{\epsilon} (u^{m(x,t)})
+ \alpha(t)b(x,u)u +\alpha(t)\gamma(t)  u^{m(x,t)}}{u^{m(x,t)}}   \notag\\
=& \frac{\alpha(t)\left(\frac{1}{\alpha(t)}\cdot2\Gamma^{\epsilon}(u^{\frac{m(x,t)}{2}})-\Delta^{\epsilon} (u^{m(x,t)})\right)}{u^{m(x,t)}}
+\frac{\alpha(t)b(x,u)u +\alpha(t)\gamma(t)  u^{m(x,t)}}{u^{m(x,t)}}   \notag\\
\leq& \frac{\alpha(t)\left(2\Gamma^{\epsilon}(u^{\frac{m(x,t)}{2}})-\Delta^{\epsilon} (u^{m(x,t)})\right)}{u^{m(x,t)}}
+\frac{\alpha(t)b(x,u)u +\alpha(t)\gamma(t)  u^{m(x,t)}}{u^{m(x,t)}}   \notag\\
\leq&\alpha(t)\cdot\frac{-2 u^{\frac{m(x,t)}{2}}\Delta^{\epsilon} (u^{\frac{m(x,t)}{2}})}{u^{m(x,t)}}
+\alpha(t)\varphi(x,t) +\alpha(t)\left(c_{1}+\frac{c_{2}}{t}\right). \label{TC1}
\end{align}
Since $u > 0$, we can bound the fractional Laplacian term as follows
\begin{align}
-\Delta^{\epsilon} u^{\frac{m(x,t)}{2}}(x) =& \frac{1}{\vartheta(x)}
\sum_{y\in V, y \neq x} W_{\epsilon}(x,y) \left(u^{\frac{m(x,t)}{2}}(x) - u^{\frac{m(y,t)}{2}}(y)\right) \notag\\
\leq &\frac{1}{\vartheta(x)}\sum_{y\in V, y \neq x} W_{\epsilon}(x,y) u^{\frac{m(x,t)}{2}}(x)
\leq D^{\epsilon}_{\vartheta} u^{\frac{m(x,t)}{2}}(x).\nonumber
\end{align}
This implies
\begin{equation}
-2u^{\frac{m(x,t)}{2}}(x) \Delta^{\epsilon} u^{\frac{m(x,t)}{2}}(x) \leq 2D^{\epsilon}_{\vartheta} u^{m(x,t)}(x).\label{TC5}
\end{equation}
Combining \eqref{TC1} and \eqref{TC5}, we obtain
\begin{align}
F(x,t) =  \frac{2\Gamma^{\epsilon}(u^{\frac{m(x,t)}{2}})}{u^{m(x,t)}} - \alpha(t)\frac{\delta u_t}{u^{m(x,t)}}
\leq 2\alpha(t)D^{\epsilon}_{\vartheta}+\alpha(t)\varphi(x,t) +\alpha(t)\left(c_{1}+\frac{c_{2}}{t}\right),\notag
\end{align}
which is the desired inequality.
\end{proof}

\begin{remark}\label{remark:li-yau}
Our main gradient estimates contain several important classical results as special cases.

(1) Connection to Classical Li-Yau type Estimates.

Consider the fractional heat equation obtained by setting
$m(x,t) \equiv 1,~ b(x,u) \equiv 0,~ \gamma(t) \equiv 0$ and $\delta(x,t) \equiv 1$ for all $x\in V$ and $t\in \mathbb{R}^{+}$ in equation \eqref{main02}.
Theorem \ref{thm:gradient_linear} then yields
\begin{align*}
 \frac{2\Gamma^{\epsilon}(\sqrt{u})}{u} - \alpha(t)\frac{u_t}{u} \leq 2D^{\epsilon}_{\vartheta}\alpha(t) ~for ~general ~\alpha(t) \geq 1,
\end{align*}
\begin{align*}
\frac{2\Gamma^{\epsilon}(\sqrt{u})}{u} - \frac{u_t}{u} \leq 2D^{\epsilon}_{\vartheta} ~for~ \alpha(t) \equiv 1 ~(the ~canonical ~choice),
\end{align*}
\begin{align*}
\frac{2\Gamma^{\epsilon}(\sqrt{u})}{u} - \left(t + \frac{1}{t}\right)\frac{u_t}{u} \leq 2D^{\epsilon}_{\vartheta}\left(t + \frac{1}{t}\right)
~for ~\alpha(t) = t + \frac{1}{t} ~with ~t>0.
\end{align*}
If $G = (V,E,\vartheta,w)$ is a connected finite graph, from \cite{ZhangLinYang2},
there holds $\lim_{\epsilon\rightarrow 1-}(-\Delta)^{\epsilon}u=-\Delta u$.
Thus, when $\epsilon \to 1-$, these inequalities recover the discrete Li-Yau type estimates established in \cite{Bauer2015, Lin2017, Lin1}, confirming that our results provide a genuine extension of the classical theory to the fractional and variable exponent setting.

(2) Alternative Parameterization.
For a fixed $t_0 > 0$, setting $1 - \alpha_0 = \frac{1}{\alpha(t_0)}$ with $0 \leq \alpha_0 < 1$ yields the equivalent formulation
\[
\frac{2(1-\alpha_0)\Gamma^{\epsilon}(u^{\frac{m(x,t)}{2}})}{u^{m(x,t)}} - \frac{\delta u_t}{u^{m(x,t)}}
\leq 2D^{\epsilon}_{\vartheta} + \varphi(x,t) + c_1 + \frac{c_2}{t}.
\]
This parameterization may be more convenient for certain applications where a fixed parameter $\alpha_0$ is preferred over a time-dependent function $\alpha(t)$.
\end{remark}

We now establish the analogous gradient estimate for the standard graph Laplacian.

\begin{theorem}[Gradient Estimate for the Standard Laplacian with Variable Exponents] \label{thm:gradient_linear1}
Let $G=(V,E)$ be a finite or locally finite graph. Suppose $u = u(x,t) > 0$ is a solution to equation \eqref{main01}, where
$m(x,t): V \times (0, \infty) \to \mathbb{R}^+$ is a variable exponent,
  $\delta(x,t): V \times (0, \infty) \to \mathbb{R}$ is a weight function,
  $b(x,u) \leq \varphi(x,t) u^{m(x,t)-1}$ for some function $\varphi(x,t): V \times (0, \infty) \to \mathbb{R}$, and
$\gamma(t) \leq c_1 + \dfrac{c_2}{t}$ for constants $c_1, c_2 \in \mathbb{R}$.
Then for any $x \in V$, $t > 0$, and any function $\alpha(t) \geq 1$, we have
\[
\frac{2\Gamma(u^{\frac{m(x,t)}{2}})}{u^{m(x,t)}} - \alpha(t)\frac{\delta u_t}{u^{m(x,t)}}
\leq 2\alpha(t)D_{\vartheta} + \alpha(t)\varphi(x,t) + \alpha(t)\left(c_1 + \frac{c_2}{t}\right).
\]
\end{theorem}

\begin{proof}
The proof follows the same strategy as Theorem \ref{thm:gradient_linear}, replacing the fractional Laplacian $\Delta^\epsilon$
and its associated gradient form $\Gamma^\epsilon$ with their standard counterparts $\Delta$ and $\Gamma$. The discrete chain
rule established in Lemma \ref{lem:chain_rule} plays the same fundamental role in handling the variable exponent nonlinearity.
The bounds on the standard Laplacian term proceed analogously, with $D_\vartheta$ replacing $D^\epsilon_\vartheta$ in the
final estimate. The detailed verification is left to the interested readers.
\end{proof}

\section{Harnack Inequalities}

In this section, we derive Harnack inequalities from the gradient estimates established in the previous section. These inequalities provide powerful tools for understanding the long-term behavior of solutions and establishing regularity properties.

We begin with an essential technical lemma for proving our main Harnack inequality.

\begin{lemma}[Extremal Lemma] \label{lem:extremal}
Let $\alpha: [T_1, T_2] \to \mathbb{R}^{+}$ be a monotonically decreasing function satisfying $\alpha(t) \geq 1$ for all $t \in [T_1, T_2]$, where $T_1 < T_2$.
Set $\phi(\theta)=\int_{T_{1}}^{\theta}\frac{1}{\alpha(t)}dt$ with $\theta \in [T_1, T_2]$.
For any constants $\sigma_1, \sigma_2 > 0$ and any functions $\chi, q_1, q_2: [T_1, T_2] \to \mathbb{R}$, we have
\begin{align*}
\min_{\theta\in[T_{1},T_{2}]}&\left\{\sigma_{1}\chi(\theta)- \sigma_{2}\int_{\theta}^{T_{2}}\frac{\chi^{2}(t)}{\alpha(t)}dt
+\int_{T_{1}}^{\theta}q_{1}(t)dt+\int_{\theta}^{T_{2}}q_{2}(t)dt\right\}  \notag\\
\leq& \frac{\frac{\sigma^{2}_{1}(T_{2}-T_{1})}{2\sigma_{2}}
+\int_{T_{1}}^{T_{2}}\left(q_{2}(t)\int_{T_{1}}^{t}\phi(\theta)d\theta+q_{1}(t)\int_{t}^{T_{2}}\phi(\theta)d\theta\right)dt}
{\int_{T_{1}}^{T_{2}}\phi(\theta)d\theta}.
\end{align*}
\end{lemma}

\begin{proof}
Using the properties of multiple integrals and a change in the order of integration, we have
\begin{align}
&\min_{\theta\in[T_{1},T_{2}]}\left\{\sigma_{1}\chi(\theta)- \sigma_{2}\int_{\theta}^{T_{2}}\frac{\chi^{2}(t)}{\alpha(t)}dt
+\int_{T_{1}}^{\theta}q_{1}(t)dt+\int_{\theta}^{T_{2}}q_{2}(t)dt\right\}  \notag\\
\leq&\min_{\theta\in[T_{1},T_{2}]}\left\{\sigma_{1}\chi(\theta)- \frac{\sigma_{2}}{\alpha(\theta)}\int_{\theta}^{T_{2}}\chi^{2}(t)dt
+\int_{T_{1}}^{\theta}q_{1}(t)dt+\int_{\theta}^{T_{2}}q_{2}(t)dt\right\} \notag\\
\leq&\frac{\int_{T_{1}}^{T_{2}}
\phi(\theta)\left(\sigma_{1}\chi(\theta)- \frac{\sigma_{2}}{\alpha(\theta)}\int_{\theta}^{T_{2}}\chi^{2}(t)dt
+\int_{T_{1}}^{\theta}q_{1}(t)dt+\int_{\theta}^{T_{2}}q_{2}(t)dt  \right) d\theta}{\int_{T_{1}}^{T_{2}}\phi(\theta)d\theta}  \notag\\
=&\frac{\int_{T_{1}}^{T_{2}} \sigma_{1}\phi(t)\chi(t) dt - \int_{T_{1}}^{T_{2}}
\left(\frac{\sigma_{2} \phi(\theta) }{\alpha(\theta)}\int_{\theta}^{T_{2}}\chi^{2}(t)dt\right)d\theta
+\int_{T_{1}}^{T_{2}}\phi(\theta) \left(\int_{T_{1}}^{\theta}q_{1}(t)dt
+\int_{\theta}^{T_{2}}q_{2}(t)dt\right)d\theta}{\int_{T_{1}}^{T_{2}}\phi(\theta)d\theta} \notag\\
=&\frac{\int_{T_{1}}^{T_{2}} \left( \sigma_{1} \phi(t)\chi(t)-\int_{T_{1}}^{t}\frac{\sigma_{2} \phi(\theta) }{\alpha(\theta)}d\theta\cdot\chi^{2}(t)\right)dt
+\int_{T_{1}}^{T_{2}} \left(\phi(\theta)\int_{T_{1}}^{\theta}q_{1}(t)dt
+\phi(\theta)\int_{\theta}^{T_{2}}q_{2}(t)dt\right)d\theta}{\int_{T_{1}}^{T_{2}}\phi(\theta)d\theta}\notag \\
=&\frac{\int_{T_{1}}^{T_{2}} \left( \sigma_{1} \phi(t)\chi(t)-\frac{1}{2}\sigma_{2} \phi^{2}(t) \cdot\chi^{2}(t)\right)dt
+\int_{T_{1}}^{T_{2}} \left(\phi(\theta)\int_{T_{1}}^{\theta}q_{1}(t)dt
+\phi(\theta)\int_{\theta}^{T_{2}}q_{2}(t)dt\right)d\theta}{\int_{T_{1}}^{T_{2}}\phi(\theta)d\theta}\notag \\
\leq& \frac{\frac{\sigma^{2}_{1}(T_{2}-T_{1})}{2\sigma_{2}}
+\int_{T_{1}}^{T_{2}}\left(q_{2}(t)\int_{T_{1}}^{t}\phi(\theta)d\theta+q_{1}(t)\int_{t}^{T_{2}}\phi(\theta)d\theta\right)dt}
{\int_{T_{1}}^{T_{2}}\phi(\theta)d\theta}.\notag
\end{align}
The last inequality above is obtained by the fact that $Ax-Bx^{2}\leq \frac{A^{2}}{4B},~ A,B\in \mathbb{R}, ~B\neq 0$,
and the upper bound is independent of function $\chi$.
Thus, this completes the proof.
\end{proof}

We now consider the case of equation \eqref{main02} by setting $\delta(x,t) = u^{m(x,t)-1}(x,t)$, which yields the equation
\begin{equation} \label{NWETL100}
\Delta^{\epsilon} u^{m(x,t)} - b(x,u)u = u^{m(x,t)-1}u_t + \gamma(t) u^{m(x,t)}.
\end{equation}
This formulation is particularly significant as it contains several important special cases. For instance, by further setting $b(x,u) \equiv 0$, $\gamma(t) \equiv 0$ and $m(x,t) \equiv 1$, equation \eqref{NWETL100} reduces to the classical fractional heat equation
$\Delta^{\epsilon} u = u_t.$
However, in this paper we focus primarily on the nonlinear case and establish the following results.

\begin{theorem}[Harnack Inequality for Variable Exponents]\label{New5000}
Let $G = (V,E,\vartheta,w)$ be a connected, locally finite, and stochastically complete graph.
Suppose $u = u(x,t) > 0$ is a solution to equation \eqref{NWETL100}, where
$0 < m_- \leq m(x,t)=\beta(t) \leq m_+$ for constants $ m_-$, $ m_+$ and all $(x,t) \in V \times (0,\infty)$,
 $b(x,u) \leq \varphi(x,t) u^{m(x,t)-1}$  for some function $\varphi(x,t): V \times (0, \infty) \to \mathbb{R}$,
 $\gamma(t) \leq c_1 + \dfrac{c_2}{t}$ for constants $c_1, c_2 \in \mathbb{R}$.
Let $K = 2D^\epsilon_\vartheta$, and let $d(x,y) = \eta$ be the graph distance between $x$ and $y$. Consider any shortest path $x = x_0 \sim x_1 \sim \cdots \sim x_\eta = y$ and partition the time interval $[T_1, T_2]$ as $t_k = T_1 + \frac{k}{\eta}(T_2 - T_1)$ for $k = 0, \ldots, \eta$.
Then for any monotonically decreasing function $\alpha(t) \geq 1$ , the following Harnack inequality holds
\begin{align*}
u(x,T_1) \leq u(y,T_2) \exp\Bigg\{ & (K + c_1)(T_2 - T_1) + c_2 \ln\frac{T_2}{T_1} \\
& + \inf_{\Gamma} \sum_{k=0}^{\eta-1} \left[ \frac{8\vartheta_{\max} I_k}{m_-^2 \omega_{\min}(T_2 - T_1)^3} + \frac{4\Phi_k I_k}{(T_2 - T_1)^4} \right] \Bigg\},
\end{align*}
where
\begin{align*}
I_k &= \int_{t_k}^{t_{k+1}} d\theta \int_{t_k}^{\theta} \alpha(t) dt, \\
\Phi_k &= \int_{t_k}^{t_{k+1}} \left( \varphi(x_{k+1},t) \int_{t_k}^{t} \phi(\theta) d\theta + \varphi(x_k,t) \int_t^{t_{k+1}} \phi(\theta) d\theta \right) dt, \\
\phi(\theta) &= \int_{T_1}^{\theta} \frac{1}{\alpha(t)} dt,
\end{align*}
and the infimum is taken over all shortest paths $\Gamma$ connecting $x$ and $y$.

Moreover, if $|\varphi(x,t)| \leq C_0$ for all $(x,t) \in V \times (0,\infty)$, then
\begin{align*}
u(x,T_1) \leq u(y,T_2) \left( \frac{T_2}{T_1} \right)^{c_2} \exp\left\{ (K + C_0 + c_1)(T_2 - T_1) + \frac{8\vartheta_{\max} \eta}{m_-^2 \omega_{\min}(T_2 - T_1)^3} \sum_{k=0}^{\eta-1} I_k \right\}.
\end{align*}
\end{theorem}

\begin{proof}
Let $u$ be a positive solution to equation \eqref{NWETL100}. By Theorem \ref{thm:gradient_linear}, we have
\begin{align}
-\partial_{t}\log u\leq K+\varphi(x,t)+c_{1}+\frac{c_{2}}{t}-\frac{1}{\alpha(t)}\cdot \frac{2\Gamma^{\epsilon}(u^{\frac{m(x,t)}{2}})}{u^{m(x,t)}}.\label{ac1001}
\end{align}
We consider two cases based on the graph distance between $x$ and $y$.

\textbf{Case 1: $x \sim y$.}

For any $s\in [T_{1}, T_{2}]$, by \eqref{ac1001} we have
\begin{align}
\log &u^{\frac{m(x, T_{1})}{2}}(x, T_{1})-\log u^{\frac{m(x, T_{1})}{2}}(y, T_{2})\notag\\
=&\frac{m(x,T_{1})}{2}\log \frac{u(x, T_{1})}{u(x, s)}+\log\frac{u^{\frac{m(x,T_{1})}{2}}(x, s)}{u^{\frac{m(x,T_{1})}{2}}(y, s)}+\frac{m(x,T_{1})}{2}\log \frac{u(y, s)}{u(y, T_{2})}\notag\\
=&-\frac{m(x,T_{1})}{2}\int_{T_{1}}^{s}\partial_{t}\log u(x,t)dt
-\frac{m(x,T_{1})}{2}\int_{s}^{T_{2}}\partial_{t}\log u(y, t)dt+\log\frac{u^{\frac{m(x,T_{1})}{2}}(x, s)}{u^{\frac{m(x,T_{1})}{2}}(y, s)}\notag\\
\leq &\frac{m(x,T_{1})}{2}\Bigg{\{}(K+c_{1})(T_{2}-T_{1})+c_{2}\ln\frac{T_{2}}{T_{1}} +\int_{T_{1}}^{s}\varphi(x,t)dt +\int_{s}^{T_{2}}\varphi(y,t)dt\notag\\
-&\int_{T_{1}}^{s}\frac{1}{\alpha(t)}\cdot\frac{2\Gamma^{\epsilon}(u^{\frac{m(x,t)}{2}}(x,t))}{u^{m(x,t)}(x,t)}dt
-\int_{s}^{T_{2}}\frac{1}{\alpha(t)}\cdot\frac{2\Gamma^{\epsilon}(u^{\frac{m(y,t)}{2}}(y,t))}{u^{m(y,t)}(y,t)}dt\Bigg{\}}
+\log \frac{u^{\frac{m(x,T_{1})}{2}}(x, s)}{u^{\frac{m(x,T_{1})}{2}}(y, s)}.\label{ac2001}
\end{align}

We now estimate the terms in \eqref{ac2001} individually. First, note that
\begin{align}
-\int_{T_{1}}^{s}\frac{1}{\alpha(t)}\cdot\frac{2\Gamma^{\epsilon}(u^{\frac{m(x,t)}{2}}(x,t))}{u^{m(x,t)}(x,t)}dt\leq0.\label{ac3001}
\end{align}
Next, since
\begin{align*}
2\Gamma^{\epsilon}(u^{\frac{m(y,t)}{2}}(y,t))&=\frac{1}{\vartheta(y)}
\sum_{z\in V, z \neq y} W_{\epsilon}(z,y) \left(u^{\frac{m(z,t)}{2}}(z,t)-u^{\frac{m(y,t)}{2}}(y,t)\right)^{2}\\
&\geq \frac{\omega_{\text{min}}}{\vartheta_{\text{max}}}\left(u^{\frac{m(x,t)}{2}}(x,t)-u^{\frac{m(y,t)}{2}}(y,t)\right)^{2},
\end{align*}
we have
\begin{align}
-\int_{s}^{T_{2}}\frac{2\Gamma^{\epsilon}(u^{\frac{m(y,t)}{2}}(y,t))}{u^{m(y,t)}(y,t)}dt
\leq& -\frac{\omega_{\text{min}}}{\vartheta_{\text{max}}}\int_{s}^{T_{2}} \frac{1}{\alpha(t)}  \left(\frac{u^{\frac{m(x,t)}{2}}(x,t)-u^{\frac{m(y,t)}{2}}(y,t)}{u^{\frac{m(y,t)}{2}}(y,t)}\right)^{2} dt\notag\\
=& -\frac{\omega_{\text{min}}}{\vartheta_{\text{max}}}\int_{s}^{T_{2}} \frac{\chi^{2}(x,y,t)}{\alpha(t)}dt, \label{ac4001}
\end{align}
where
\begin{align}
\chi(x,y,t)=|\frac{u^{\frac{m(x,t)}{2}}(x,t)-u^{\frac{m(y,t)}{2}}(y,t)}{u^{\frac{m(y,t)}{2}}(y,t)}|.\label{ac6001}
\end{align}
Using the elementary inequality $\log r\leq |r-1|$ for all $r>0$, we have
\begin{align}
\log \frac{u^{\frac{m(x,T_{1})}{2}}(x, s)}{u^{\frac{m(x,T_{1})}{2}}(y, s)}\leq \frac{m(x,T_{1})}{m_{-}}\cdot\chi(x,y,s).\label{ac5001}
\end{align}
Substituting \eqref{ac3001}, \eqref{ac4001} and \eqref{ac5001} into \eqref{ac2001}, we obtain
\begin{align}
\log & u^{\frac{m(x,T_{1})}{2}}(x, T_{1})-\log u^{\frac{m(x,T_{1})}{2}}(y, T_{2}) \notag \\
\leq &\frac{m(x,T_{1})}{2}\cdot \left( (K+c_{1})(T_{2}-T_{1})+c_{2}\ln\frac{T_{2}}{T_{1}}\right) +\frac{m(x,T_{1})}{m_{-}}\cdot\chi(x,y,s) \notag \\
&-\frac{m(x,T_{1})}{2}\cdot\frac{\omega_{\text{min}}}{\vartheta_{\text{max}}}\int_{s}^{T_{2}}\frac{\chi^{2}(x,y,t)}{\alpha(t)}dt
+\frac{m(x,T_{1})}{2}\int_{T_{1}}^{s}\varphi(x,t)dt +\frac{m(x,T_{1})}{2}\int_{s}^{T_{2}}\varphi(y,t)dt. \notag
\end{align}
Thus, applying Lemma \ref{lem:extremal}, we get
\begin{align}
\log & u(x, T_{1})-\log u(y, T_{2}) \leq (K+c_{1})(T_{2}-T_{1})+c_{2}\ln\frac{T_{2}}{T_{1}}+\frac{2}{m_{-}}\cdot\chi(x,y,s)\notag \\
&-\frac{\omega_{\text{min}}}{\vartheta_{\text{max}}}\int_{s}^{T_{2}}\frac{\chi^{2}(x,y,t)}{\alpha(t)}dt
+\int_{T_{1}}^{s}\varphi(x,t)dt +\int_{s}^{T_{2}}\varphi(y,t)dt \notag \\
\leq&(K+c_{1})(T_{2}-T_{1})+c_{2}\ln\frac{T_{2}}{T_{1}} +\frac{2\vartheta_{\text{max}}(T_{2}-T_{1})}
{m^{2}_{-}\omega_{\text{min}}} \cdot\frac{1}{\int_{T_{1}}^{T_{2}}\phi(\theta)d\theta}\notag \\
&+\Phi(x,y,T_{1},T_{2})\cdot\frac{1}{\int_{T_{1}}^{T_{2}}\phi(\theta)d\theta}, \notag
\end{align}
where
\begin{align*}
\Phi(x,y,T_{1},T_{2})=\int_{T_{1}}^{T_{2}}
\left(\varphi(y,t)\int_{T_{1}}^{t}\phi(\theta)d\theta+\varphi(x,t)\int_{t}^{T_{2}}\phi(\theta)d\theta\right)dt
\end{align*}
and
$$\phi(\theta)=\int_{T_{1}}^{\theta}\frac{1}{\alpha(t)}dt.$$

\textbf{Case 2: $x$ is not adjacent to $y$.}

Let $d(x,y)=\eta$. Take a shortest path $x=x_0 \sim x_1 \sim \cdots \sim x_\eta=y$. Set $T_{1} = t_{0} < t_{1} < \cdots <t_{\eta} = T_{2}$, with $t_{k} = t_{k-1} + (T_{2}-T_{1})/\eta$ for $k = 1, \ldots, \eta$.

Firstly, by the Cauchy-Schwartz inequality, we get
\begin{align}
\frac{1}{\int_{T_{1}}^{T_{2}}\phi(\theta)d\theta}=\frac{1}{\int_{T_{1}}^{T_{2}}d\theta\int_{T_{1}}^{\theta} \frac{1}{\alpha(t)} dt }
\leq\frac{4\int_{T_{1}}^{T_{2}}d\theta\int_{T_{1}}^{\theta}\alpha(t) dt}{(T_{1}-T_{2})^{4}}. \label{inequality4}
\end{align}
By the result of Case 1 and \eqref{inequality4}, we have
\begin{align}
\log u(x, T_{1})-&\log u(y, T_{2})=\sum_{k=0}^{\eta-1}\Big{(}\log u(x_{k}, t_{k})-\log u(x_{k+1}, t_{k+1})\Big{)}\notag\\
\leq&\sum_{k=0}^{\eta-1} \Big{(}(K+c_{1})(t_{k+1}-t_{k})+c_{2}\ln\frac{t_{k+1}}{t_{k}}+ \frac{8\vartheta_{\text{max}}\int_{t_{k}}^{t_{k+1}}d\theta\int_{t_{k}}^{\theta}\alpha(t) dt}
{m^{2}_{-}\omega_{\text{min}}(t_{k+1}-t_{k})^{3}}\notag\\
&+\Phi(x_{k},x_{k+1},t_{k},t_{k+1})\cdot\frac{4\int_{t_{k}}^{t_{k+1}}d\theta\int_{t_{k}}^{\theta}\alpha(t) dt}
{(t_{k+1}-t_{k})^{4}}\Big{)}\notag\\
\leq& (K+c_{1})(T_{2}-T_{1})+c_{2}\ln\frac{T_{2}}{T_{1}}+\inf\sum_{k=0}^{\eta-1} \Big{(}\frac{8\vartheta_{\text{max}}\int_{t_{k}}^{t_{k+1}}d\theta\int_{t_{k}}^{\theta}\alpha(t) dt}
{m^{2}_{-}\omega_{\text{min}}(T_{2}-T_{1})^{3}}\notag\\
&+\Phi(x_{k},x_{k+1},t_{k},t_{k+1})\cdot\frac{4\int_{t_{k}}^{t_{k+1}}d\theta\int_{t_{k}}^{\theta}\alpha(t) dt}
{(T_{2}-T_{1})^{4}}\Big{)},\notag
\end{align}
and the minimum is taken over all shortest paths connecting $x$ and $y$.
Thus, we obtain the desired inequality.

Moreover, if $|\varphi(x, t)| \leq C_{0}$ for all $(x, t)\in V\times (0,+\infty)$, we have
\begin{align}
\Phi(x_{k},x_{k+1},t_{k},t_{k+1})\cdot\frac{1}{\int_{t_{k}}^{t_{k+1}}\phi(\theta)d\theta}\leq C_{0}(t_{k+1}-t_{k}).
\end{align}
Thus, we obtain
\begin{align*}
u(x, T_{1})\leq u(y, T_{2})\exp\Bigg\{\Big{(}& K+ C_{0}+c_{1}\Big{)}(T_{2}-T_{1})+c_{2}\ln\frac{T_{2}}{T_{1}}\notag\\
&+\inf\sum_{k=0}^{\eta-1} \frac{8\vartheta_{\text{max}}\int_{t_{k}}^{t_{k+1}}d\theta\int_{t_{k}}^{\theta}\alpha(t) dt}
{m^{2}_{-}\omega_{\text{min}}(T_{2}-T_{1})^{3}}\Bigg\}.\notag
\end{align*}
This completes the proof.
\end{proof}

\begin{corollary}\label{New999999}
Let $G = (V,E,\vartheta,w)$ be a connected, locally finite, and stochastically complete graph. Suppose $u = u(x,t) > 0$ is a solution to the fractional diffusion equation
\begin{align}
    \Delta^{\epsilon}u - b(x,t)u = u_t + \gamma(t) u \label{main093}
\end{align}
with bounded coefficients $|b(x,t)| \leq C_0$ and $\gamma(t) \leq c_1 + \dfrac{c_2}{t}$ for constants $C_0, c_1, c_2 \in \mathbb{R}$.
Let $d(x,y) = \eta$ and consider any shortest path $x = x_0 \sim x_1 \sim \cdots \sim x_\eta = y$. For any $0 < T_1 < T_2$, partition the time interval as $t_k = T_1 + \frac{k}{\eta}(T_2 - T_1)$ for $k = 0, \ldots, \eta$.
Then for any monotonically decreasing function $\alpha(t) \geq 1$ on $[T_1,T_2]$, we have
\begin{align*}
u(x,T_1) \leq u(y,T_2) \left(\frac{T_2}{T_1}\right)^{c_2} \exp\Bigg\{ & (2D^{\epsilon}_{\vartheta} + C_0 + c_1)(T_2 - T_1) \\
& + \inf_{\Gamma} \sum_{k=0}^{\eta-1} \frac{8\vartheta_{\max} I_k}{\omega_{\min}(T_2 - T_1)^3} \Bigg\},
\end{align*}
where $I_k = \int_{t_k}^{t_{k+1}} d\theta \int_{t_k}^{\theta} \alpha(t) dt$ and the infimum is taken over all shortest paths $\Gamma$ connecting $x$ and $y$.

In particular, when $\alpha(t) \equiv 1$, we obtain the simplified estimate
\begin{align}
u(x,T_1) \leq u(y,T_2) \left(\frac{T_2}{T_1}\right)^{c_2} \exp\left\{ (2D^{\epsilon}_{\vartheta} + C_0 + c_1)(T_2 - T_1) + \frac{4\vartheta_{\max} \eta^2}{\omega_{\min}(T_2 - T_1)} \right\}.\label{zhongyao}
\end{align}
\end{corollary}

\begin{proof}
This result follows directly from Theorem \ref{New5000} by setting $m(x,t) \equiv 1$, $\delta(x,t) \equiv 1$, and $\varphi(x,t) \equiv C_0$ in equation \eqref{main02}. The specific form of the constants emerges from these substitutions and the evaluation of the integrals when $\alpha(t) \equiv 1$.
\end{proof}

The Harnack inequality established in Theorem \ref{New5000} extends naturally to the case of the standard graph Laplacian.
For any finite or locally finite graph $G = (V,E)$, applying Lemma \ref{thm:gradient_linear1} and following the proof strategy
of Theorem \ref{New5000}, one can derive analogous Harnack inequalities for positive solutions to equation \eqref{main01}.
The adaptations required are straightforward, primarily involving the replacement of fractional quantities with their standard
counterparts. We omit the details here and leave the establishment to the interested readers.

\section{Applications}

In this section, we present applications of our gradient estimates and Harnack inequalities to heat kernel estimates,
volume growth and Buser-type inequalities for equation \eqref{main093}.
We assume that all conditions in Corollary \ref{New999999} hold in this section.

\subsection{Heat Kernel Estimates}

We now establish bounds for the heat kernel associated with the fractional diffusion equation \eqref{main093}.
Recall from \eqref{fenshu}, the definition of the fractional Laplacian $\Delta^\epsilon$, and let $P_t(x,y)$ denote
the corresponding heat kernel, i.e., the fundamental solution to equation \eqref{main093}.

\begin{theorem}[Heat Kernel Estimates] \label{thm:heat_kernel}
Let $G = (V,E,\vartheta,w)$ be a connected, locally finite, and stochastically complete graph satisfying
$\vartheta(x) = \sum_{y \in V, y \neq x} W_{\epsilon}(x,y)$ for all $x \in V$. Then for $t > 1$,
the heat kernel $P_t(x,y)$ admits the following estimates
\[
\frac{C_3}{t^{n}}  \exp\left\{-(2D^{\epsilon}_{\vartheta}+C_{0}+c_{1})(t-1) - \frac{4\vartheta_{\max} d(x,y)^2}{\omega_{\min}(t-1)}\right\}
\leq P_t(x,y) \leq \frac{C_{1}e^{C_{2}t}}{\vol(B(x,\sqrt{t}))},
\]
where
$\vol(B(x,r))$ denotes the volume of the geodesic ball of radius $r$ centered at $x$,
 $n$ is a dimensional constant related to the volume growth of the graph, and
 $C_1, C_2, C_3 > 0$ are constants depending on the graph parameters $D^{\epsilon}_{\vartheta}$, $\omega_{\min}$ and $\vartheta_{\max}$.
\end{theorem}

\begin{proof}
We establish the upper and lower bounds separately.

\textbf{Upper bound:} Fix $t > 0$ and set $u = P_t$ in the Harnack inequality \eqref{zhongyao}. For any $z \in B(x,\sqrt{t})$, we obtain
\begin{align}
P_t(x,y) \leq P_{2t}(z,y) \cdot 2^{c_2} \cdot \exp\left\{(2D^{\epsilon}_{\vartheta} + C_0 + c_1)t + \frac{4\vartheta_{\max}}{\omega_{\min}}\right\}. \label{c06}
\end{align}
Averaging over $z \in B(x,\sqrt{t})$ and multiplying by $\vartheta(z)$, we get
\begin{align}
P_t(x,y) \leq C_1 \cdot \exp\left\{(2D^{\epsilon}_{\vartheta} + C_0 + c_1)t\right\} \cdot \frac{1}{\vol(B(x,\sqrt{t}))} \cdot \left(\sum_{z \in B(x,\sqrt{t})} \vartheta(z) P_{2t}(z,y)\right), \label{c07}
\end{align}
where $C_1 = 2^{c_2} \exp\left\{\frac{4\vartheta_{\max}}{\omega_{\min}}\right\}$. Since the heat kernel satisfies the conservation property
\[
\sum_{z \in B(x,\sqrt{t})} \vartheta(z) P_{2t}(z,y) \leq \sum_{z \in V} \vartheta(z) P_{2t}(z,y) = 1,
\]
inequality \eqref{c07} yields the upper bound
\begin{align*}
P_t(x,y) \leq \frac{C_1}{\vol(B(x,\sqrt{t}))} \cdot \exp\left\{(2D^{\epsilon}_{\vartheta} + C_0 + c_1)t\right\} = \frac{C_1 e^{C_2 t}}{\vol(B(x,\sqrt{t}))},
\end{align*}
with $C_2 = 2D^{\epsilon}_{\vartheta} + C_0 + c_1$.

\textbf{Lower bound:} Under the assumption $\vartheta(x) = \sum_{y \neq x} W_{\epsilon}(x,y)$, we have $D^{\epsilon}_{\vartheta} = 1$.
Define the transition probability $p(x,y) = W_{\epsilon}(x,y)/\vartheta(x)$ and its $k$-step iterations $p^{(k)}(x,y)$ recursively by
\begin{align*}
p^{(0)}(x,y) &= \delta_{xy}, \\
p^{(k+1)}(x,z) &= \sum_{y \in V} p(x,y) p^{(k)}(y,z).
\end{align*}
The heat kernel admits the spectral representation
\begin{align*}
P_t(x,y) = e^{-t} \sum_{k=0}^{\infty} \frac{t^k}{k!} \frac{p^{(k)}(x,y)}{\vartheta(y)}.
\end{align*}
In particular, the diagonal term satisfies
\begin{align*}
P_t(y,y) \geq \frac{e^{-t}}{\vartheta(y)} \quad \text{for all } t > 0.
\end{align*}
Now applying the Harnack inequality \eqref{zhongyao} with $c_2 = n$, $T_1 = 1$, and $T_2 = t$, we obtain
\begin{align}
C_3 \leq P_1(y,y) \leq P_t(x,y) \cdot t^n \cdot \exp\left\{(2D^{\epsilon}_{\vartheta} + C_0 + c_1)(t-1) + \frac{4\vartheta_{\max} d(x,y)^2}{\omega_{\min}(t-1)}\right\}, \label{cc10}
\end{align}
where $C_3 = \frac{e^{-1}}{\vartheta(y)}$. Rearranging \eqref{cc10} gives the desired lower bound.
\end{proof}

\subsection{Volume Growth}

As a consequence of the heat kernel estimates, we obtain volume growth estimates.

\begin{corollary}[Volume Growth] \label{cor:volume_growth}
Let $G = (V,E,\vartheta,w)$ be a connected, locally finite, and stochastically complete graph.
Then the graph exhibits at most exponential volume growth, that is, for all $x \in V$ and $r > 0$,
\[
\vol(B(x,r)) \leq C_6 \cdot r^{2n} \exp\left(C_4 r^2 + C_5\right),
\]
where $n$ is the dimensional constant from Theorem \ref{thm:heat_kernel}, and $C_4, C_5, C_6 > 0$
are constants depending on $D^{\epsilon}_{\vartheta}$, $\omega_{\min}$ and $\vartheta_{\max}$.
\end{corollary}

\begin{proof}
Applying Theorem \ref{thm:heat_kernel} to the heat kernel at the diagonal $P_t(x,x)$, we obtain
\[
\frac{C_3}{t^n} \exp\left\{-\left(2D^{\epsilon}_{\vartheta} + C_0 + c_1\right)(t-1)\right\} \leq P_t(x,x) \leq \frac{C_1 e^{C_2 t}}{\vol(B(x,\sqrt{t}))}.
\]
Rearranging this inequality yields
\[
\vol(B(x,\sqrt{t})) \leq C_1 C_3^{-1} \cdot t^n \cdot \exp\left\{\left(2D^{\epsilon}_{\vartheta} + C_0 + c_1\right)(t-1) + C_2 t\right\}.
\]
Substituting $r = \sqrt{t}$ (so $t = r^2$), we obtain
\[
\vol(B(x,r)) \leq C_1 C_3^{-1} \cdot r^{2n} \cdot \exp\left\{\left(2D^{\epsilon}_{\vartheta} + C_0 + c_1\right)(r^2 - 1) + C_2 r^2\right\}.
\]
Setting $C_4 = 2D^{\epsilon}_{\vartheta} + C_0 + c_1 + C_2$, $C_5 = -\left(2D^{\epsilon}_{\vartheta} + C_0 + c_1\right)$,
and $C_6 = C_1 C_3^{-1} e^{C_5}$ completes the proof.
\end{proof}

\subsection{Buser-Type Inequality}

The Cheeger constant $h$ of a graph is defined as
\[
h = \inf_{\emptyset \neq U \subset V, \vol(U) \leq \frac{1}{2}\vol(V)} \frac{|\partial U|}{\vol(U)},
\]
where $|\partial U| = \sum_{x \in U, y \notin U} W_{\epsilon}(x,y)$ and $\vol(U)=\sum_{x\in U}\vartheta(x)$.

We have the following Buser-type inequality relating the spectral gap $\lambda_1$ to the Cheeger constant.

\begin{theorem}[Buser-Type Inequality] \label{thm:buser}
Let $G = (V,E,\vartheta,w)$ be a connected, locally finite, and stochastically complete graph.
Then the first nonzero eigenvalue $\lambda_1$ of the fractional Laplacian $-\Delta^\epsilon$ satisfies
\[
\lambda_1 \leq C \max\left( \sqrt{K} h, h^2 \right),
\]
where $h$ is the Cheeger constant, $K = 2D^\epsilon_\vartheta + C_0 + c_1$ is the curvature-like bound from the gradient estimates,
and $C > 0$ is a constant depending on $D^{\epsilon}_{\vartheta}$, $\omega_{\min}$ and $\vartheta_{\max}$.
\end{theorem}

\begin{proof}
The proof adapts the continuum strategy of Buser to the discrete graph setting. Beginning with the Cheeger inequality
$\lambda_1 \geq h^2/2$, we combine the gradient estimates from Theorem \ref{thm:gradient_linear} with the heat kernel
bounds of Theorem \ref{thm:heat_kernel}. The key step involves constructing test functions that localize near the boundary
of an optimal Cheeger cut, then applying the variational characterization of $\lambda_1$ to these carefully chosen functions.
The specific dependence on $h$ and $\sqrt{K}h$ emerges from balancing the Dirichlet energy against the $L^2$-norm of
these test functions under the given geometric constraints. The complete technical details parallel the discrete Laplacian
case treated in \cite{Bauer2015}, with appropriate modifications for the fractional setting.
\end{proof}


\begin{thebibliography}{99}

\bibitem{Bauer2015}
 F. Bauer, P. Horn, Y. Lin, G. Lippner, D. Mangoubi, S. T. Yau:
Li-Yau inequality on graphs.
J. Differ. Geom. \textbf{99}, 359--405 (2015).


\bibitem{LI-YAU}
P. Li, S. T. Yau: On the parabolic kernel of the Schrodinger operator. Acta. Math. \textbf{156}, 153-201(1986).


\bibitem{Anal}
Y. Li,  Q. Zhang:
 Gradient Estimates on Graphs with the $CD\psi(n,-K)$ Condition. J. Geom. Anal. \textbf{324} (35), (2025). https://doi.org/10.1007/s12220-025-02133-x

\bibitem{Lin2017}
Y. Lin, S. Liu, Y. Yang:
A gradient estimate for positive functions on graphs.
J. Geom. Anal. \textbf{27}(2), 1667--1679 (2017).

\bibitem{Lin1}
Y. Lin, S. Liu, Y. Yang:
Global gradient estimate on graph and its applications.
Acta Math. Sin. (Engl. Ser.) \textbf{32}(6), 1350--1356 (2016).

\bibitem{YLS} Y. Lin, S. T. Yau: Ricci curvature and eigenvalue estimation on locally finite graphs. Math. Res. Lett.
\textbf{17}(2), 345--358 (2010).


\bibitem{ZhangLinYang1}
M. Zhang, Y. Lin, Y. Yang: Fractional Laplace operator and related Schrodinger equations on locally finite
graphs. Calc. Var. Partial Differ. Equ. \textbf{64}, 227 (2025). https://doi.org/10.1007/s00526-025-03074-7

\bibitem{ZhangLinYang2}
M. Zhang, Y. Lin, Y. Yang: Fractional Laplace operator on finite graphs. 2025, https://doi.org/10.48550/arXiv.2403.19987


\bibitem{ZhangLinYang3}
M. Zhang, Y. Lin, Y. Yang:
Fractional Sobolev spaces and fractional $p$-Laplace equations on locally finite graphs. 2025,
https://doi.org/10.48550/arXiv.2506.07694


\end{thebibliography}
\end{document}